# Decomposition-based Energy-based Dual-Phase Dynamics Identification for Nonlinear MDOF Systems

Sayantan Ghosh[1], Cristian López[2], Aryan Singh[1], Keegan J. Moore[1*]

[1]Daniel Guggenheim School of Aerospace Engineering, Georgia Institute of Technology, Atlanta, GA 30332
[2]Department of Applied Mathematics, University of Colorado Boulder, Boulder, CO 80309

*Corresponding author: (K.J. Moore)
E-mail address: kmoore@gatech.edu

## Abstract

System identification is an important step in modeling and evaluating vibrating structures, but many nonlinear system identification methods rely heavily on data-driven approaches that may not preserve physical consistency. This research extends the Energy-based Dual-phase Dynamics Identification (EDDI) method to multiple-degree-of-freedom (MDOF) mechanical structures undergoing nonlinear vibrations. The original EDDI framework was designed for single-degree-of-freedom (SDOF) systems and operates in two phases: the first identifies a model for internal nonconservative force, and the second captures internal conservative force. However, EDDI assumes that the potential energy is zero whenever the displacement is zero. For MDOF systems, this assumption requires all degrees of freedom (DOFs) to achieve zero displacement simultaneously, which simply occurs too infrequently in multimodal responses for direct application. To overcome this limitation, this work introduces Decomposition-based EDDI, which applies EDDI to decomposed response components to enable nonlinear system identification of MDOF systems. Wavelet-Bounded Empirical Mode Decomposition is used to extract nearly orthogonal, monochromatic intrinsic mode functions (IMFs) from the measured displacements. Each IMF is then treated as an individual SDOF oscillator and processed using EDDI to estimate its nonconservative and conservative internal forces. The IMF forces are then summed to reconstruct the total nonconservative and conservative forces acting on each physical DOF, which are subsequently used to identify the damping and stiffness models, respectively. The proposed EDDI framework is experimentally validated on a two-story tower structure with strong stiffness nonlinearity coupling the two floors. The results demonstrate the efficacy of EDDI in isolating and identifying complex, multi-modal nonlinear structural dynamics.

## Nomenclature

| | |
|---|---|
| CWT | Continuous wavelet transform |
| DEDDI | Decomposition-based Energy-based Dual-phase Dynamics Identification |
| DOF | Degree-of-freedom |
| EDDI | Energy-based Dual-phase Dynamics Identification |
| EMD | Empirical Mode Decomposition |
| IMF | Intrinsic mode function |
| LO | Linear Oscillator |
| MDOF | Multiple-degree-of-freedom |
| MSE | Mean squared error |
| NO | Nonlinear Oscillator |
| NSI | Nonlinear system identification |
| SDOF | Single-degree-of-freedom |
| SINDy | Sparse Identification of Nonlinear Dynamics |
| SIVA | System Identification via Validation and Adaptation |
| WBEMD | Wavelet-Bounded Empirical Mode Decomposition |
| $A_l$ | The $l$th function used to model the conservative (elastic) force |
| $B_n$ | The internal nonconservative force acting on the $n$th DOF |
| $B_{ni}$ | The dissipative force acting on the $i$th IMF for the $n$th DOF |
| $C_{ni}^{(j)}$ | The $j$th function used to model the dissipative force of $i$th IMF for the $n$th DOF |
| $I$ | The total number of IMFs extracted from an oscillatory signal |
| $I_n$ | The total number of IMFs extracted from the displacement time series of the $n$th DOF using WBEMD |
| $J$ | The total number of functions used to model the nonconservative (dissipative) force |
| $K_n$ | The internal conservative force acting on the $n$th DOF |
| $L$ | The total number of functions used to model the elastic force |
| $N$ | The total number of DOFs |
| $P$ | The total number of zero crossings |
| $\mathbf{Q}$ | Matrix of unscaled dissipated energies computed using the proposed damping model functions |
| $\mathbf{Q}_{ni}$ | Matrix of unscaled dissipated energies computed using the proposed model functions for the $i$th IMF of the $n$th DOF |
| $Q_{ni}^{(p,j)}$ | The scalar entry of the $p$th row and $j$th column of the matrix $\mathbf{Q}_{ni}$, representing the unscaled dissipated energy computed using the $j$th model function at the time instant $\gamma_p$ |
| $T$ | Kinetic energy |
| $T_{ni}$ | Kinetic energy of the $i$th IMF for the $n$th DOF |
| $a_l$ | The $l$th coefficient used to model the conservative (elastic) force |
| $b_g$ | The coefficient of the damper connecting the LO to the ground |
| $b_c$ | The coefficient of the damper coupling the LO and NO |
| $\mathbf{c}$ | Vector of unknown coefficients for the damping model for a generic SDOF system |
| $\mathbf{c}_{ni}$ | Column vector of unknown damping model coefficients for the $i$th IMF of the $n$th DOF, given by $\mathbf{c}_{ni} = \left[c_{ni}^{(1)}, c_{ni}^{(2)}, \dots, c_{ni}^{(J)}\right]^T$. |
| $c_{ni}^{(j)}$ | The $j$th scalar entry of the coefficient vector $\mathbf{c}_{ni}$, corresponding to the weight of the $j$th damping model function $C_{ni}^{(j)}$ |
| $e_{I+1}(t)$ | The remainder signal produced by the standard EMD |
| $i$ | Index indicating the IMF number sorted from lowest to highest frequency |
| $j$ | Index used to count the functions used to model the nonconservative (dissipative) force |
| $k_g$ | The stiffness connecting the LO to the ground |
| $k_1$ | The linear stiffness coupling the LO and NO |
| $k_3$ | The nonlinear stiffness coupling the LO and NO |
| $l$ | Index used to count the functions used to model the conservative force |
| $m_n$ | The mass of the $n$th DOF |
| $n$ | Index indicating the DOF |
| $p$ | Index used to count the number of zero crossings in a displacement signal or IMF |
| $\mathbf{r}$ | Vector of dissipated energies estimated from kinetic energy |
| $\mathbf{r}_{ni}$ | Vector of dissipated energies estimated from the kinetic energy of the $i$th IMF of the $n$th DOF, given by $\mathbf{r}_{ni} = \left[r_{ni}^{(1)}, r_{ni}^{(2)}, \dots, r_{ni}^{(P)}\right]$. |
| $r_{ni}^{(p)}$ | The estimated dissipated energy of the $i$th IMF for the $n$th DOF at the time instant $\gamma_p$ |
| $s(t)$ | Masking signal used to improve the decomposition of EMD |
| $t$ | Time |
| $u_+(t)$ | Signal composed of the current remainder signal plus the masking signal |
| $u_-(t)$ | Signal composed of the current remainder signal minus the masking signal |
| $x_n(t)$ | Displacement of the $n$th DOF |
| $x_{ni}$ | The $i$th IMF extracted from the displacement of the $n$th DOF |
| $\dot{x}_{ni}$ | The velocity of the $i$th IMF for the $n$th DOF |
| $\ddot{x}_{ni}$ | The acceleration of the $i$th IMF for the $n$th DOF |
| $y(t)$ | An oscillatory signal |
| $y_i(t)$ | The $i$th IMF extracted from $y(t)$ |
| $z_i(t)$ | The current remainder signal |
| $\alpha$ | Parameter used in WBEMD to scale the amplitude of the masking signal |
| $\beta$ | Parameter used in WBEMD to scale the frequency of the masking signal |
| $\gamma_p$ | The time instants where an IMF is equal to zero |
| $\omega_i$ | The characteristic frequency of the $i$th harmonic component and the $i$th IMF |

## 1. Introduction

Vibrating systems are common across many engineering domains, and their behavior can be complex especially if they exhibit nonlinearities and are coupled. These nonlinearities complicate design and analysis, requiring a deep understanding of their effects on the dynamics of the system for accurate modeling and identification. System identification is used across engineering disciplines to predict system behavior and serves as a critical tool for producing mathematical models for complex systems. Consequently, researchers have actively pursued nonlinear system identification, as reviewed in the articles [1,2], while linear system methods, such as experimental modal analysis, are already well established within structural dynamics and vibration [3]. In nonlinear structures, the behavior manifests as non-stationary transience, such as amplitude-dependent frequency and damping. Because of this, linear methods based on Fourier transformations cannot be directly applied, and the presence of nonlinear coupling in multi-degree of freedom (MDOF) systems make these traditional approaches nearly futile. To address this challenge, this work provides a framework for the nonlinear identification of MDOF systems that contain unknown nonlinearities either in stiffness or damping and are coupled via nonlinear structural elements.

Recent advancements in nonlinear system identification have increasingly leveraged data-driven methods. These techniques learn system dynamics directly from measurements using optimization, curve fitting, and combinations of nonlinear candidate functions. Typically, these methods use a large library of candidate functions, identify the unknown parameters, and eliminate the negligible terms using thresholding criteria. Although recent trends in nonlinear system identification are shifting towards data driven methods, a short overview here provides analysts with the context needed to choose the right tools. System identification methods can be broadly categorized into three types:

1. Parametric methods are preferred when the underlying physics and mathematical framework of the system are already known, reducing the identification to estimate the unknown model parameters. Examples include auto-regressive models, auto-regressive moving average, nonlinear auto-regressive moving average [4], Bayesian estimation methods [5,6], machine learning methods [7–13] and time series optimization [14].
2. Non-parametric methods are used when the physics of the problem is unknown and must be identified from the measured data directly. Examples include restoring force method [15,16] and frequency-energy response methods [17,18].
3. Semi-parametric models combine the strengths of both parametric and non-parametric techniques, and many modern data-driven techniques fall into this category. These include methods such as the Characteristic Nonlinear System Identification [19–21], Bayesian methods [22], and the widely used Sparse Identification of Nonlinear

Dynamics (SINDy) [23–26].

A more recent addition is the Energy-based Dual-phase dynamics Identification (EDDI) [27–29], which is explicitly built around energy principles that ensure compliance with the physical balance laws pertaining to an oscillating system. Additionally, an emerging class of techniques leverages generative AI, such as the System Identification via Validation and Adaptation (SIVA)[13,30]. SIVA uses the same idea as generative adversarial networks by pitting two networks, a parameter generator and validator, in competition with each other to identify mathematical models using measured response. Although the recent trends are pushing towards a heavier focus on data, energy-based methods such as EDDI and its variants are advantageous because they embed the fundamental physics into the identification process through compliant balance laws. Consequently, these methods have been shown to outperform techniques that don't leverage the inherent physics, such as SINDy [23], for the measured response of strongly nonlinear mechanical systems undergoing vibrations. They also are more resistant to measurement noise and ensure that the identified models generalize well to unseen datasets [29].

Conversely, the application of classical modal decomposition is severely limited for nonlinear MDOF systems, even in the presence of weak nonlinearities. This limitation stems from the complex phenomena introduced by nonlinearities, such as inter-mode energy interactions, internal resonances, and mode localization [31], thereby necessitating additional intermediate steps before applying system identification. To process the complex spectra of multimodal responses, numerous advanced signal processing and decomposition methods have been developed, including Empirical Mode Decomposition (EMD) [32–35], Variational Mode Decomposition [36], and Dynamic Mode Decomposition wavelet transforms [37–40]. These techniques find wide application in structural health monitoring of bolted structures [41,42], fault diagnosis in rotating machinery [43], and reduced-order modeling of the nonlinear modal interactions [44–46].

Previous applications of EDDI successfully demonstrated its efficacy on SDOF nonlinear systems using free response measurements [27], systems containing non-smooth nonlinearities [28], forced response measurements [29]. This work extends the original EDDI framework to nonlinearly coupled MDOF systems through response decomposition. Specifically, the Wavelet-Bounded Empirical Mode Decomposition (WBEMD) [34,35] to extract nearly orthogonal, monochromatic intrinsic mode functions (IMFs) from the measured displacement responses. Each IMF is then treated as the response of an SDOF oscillator, allowing EDDI to be applied directly to each component. We term this extension Decomposition-based EDDI (DEDDI). The proposed DEDDI framework is demonstrated on an experimental two-story tower structure featuring nonlinear coupling between the floors.

The remainder of this paper is structured as follows. Section 2 establishes the theoretical basis of the framework, presenting the core EDDI method, the WBEMD decomposition pipeline, and the mathematical integration required for general MDOF systems with unknown forms of nonlinearities. Section 3 details the application of the framework to the experimental setup, where the identified parameters are used to reconstruct the free response of the two-story tower under various impact force levels. To demonstrate its efficacy, the EDDI results are rigorously compared against several state-of-the-art methods, including SINDy, time-series optimization, and SIVA. Finally, Section 4 presents concluding remarks, discusses the robustness and limitations of the proposed method, and outlines potential avenues for future research.

## 2. Background & The Proposed Method

In this section, we introduce and briefly discuss EDDI and WBEMD, which are the two core components used to form DEDDI. We then explain the DEDDI in full detail. Due to the large number of indices, we will use subscripts for the indices used in WBEMD and superscripts surrounded by parentheses to indicate the indices used in EDDI.

### 2.1. Energy-based Dual-phase Dynamics Identification (EDDI)

EDDI is a two-phase nonlinear system identification method that leverages conservation of energy to identify the equation of motion for a nonlinear SDOF oscillator from measured response data. Phase I focuses on identifying a model for the nonconservative (dissipative) internal force acting on the oscillator, while the second phase produces a model for the conservative (elastic) internal force. We explain EDDI using the following equation of motion for a generic nonlinear SDOF oscillator

$$m\ddot{x} + B(x, \dot{x}) + K(x) = F(t), \tag{1}$$

where $m$ is the mass of the oscillator, $B(x, \dot{x})$ is the internal nonconservative force, $K(x)$ is the internal conservative force, $x(t)$ is the displacement, and $F(t)$ is a generic external force. Both $B(x, \dot{x})$ and $K(x)$ are assumed to be nonlinear with unknown functional form.

Phase I leverages conservation of energy to compute the energy dissipated by the internal non-conservative forces using the kinetic energy at the time instants when the potential energy is zero. To find these time instants, EDDI assumes that the potential energy is zero whenever the displacement of the oscillator is zero, then finds the time instants where the displacement is equal to zero using peak finding and linear interpolation. Let $\gamma^{(p)}$, $p = 0,1, \dots, P$ be the time instants where the displacement and potential energy are zero, then shifting time to start at $\gamma_0$, the energy dissipated can be expressed as

$$r_p = E_{diss}\left(\gamma^{(p)}\right) = T\left(\gamma^{(0)}\right) - T\left(\gamma^{(p)}\right) + W\left(\gamma^{(p)}\right), \tag{2}$$

which is expressed as

$$\int_{\gamma^{(0)}}^{\gamma^{(p)}} \dot{x}(t) B\left(x(t), \dot{x}(t)\right) dt = \frac{1}{2} m \left(\dot{x}^2\left(\gamma^{(0)}\right) - \dot{x}^2\left(\gamma^{(p)}\right)\right) + \int_{\gamma^{(0)}}^{\gamma^{(p)}} \dot{x} F(t)\, dt. \tag{3}$$

The righthand side of Eq. 3 is computed from the measured response and measured force directly, while the lefthand side is constructed by proposing a model for $B(x, \dot{x})$ of the form

$$B(x, \dot{x}) = \sum_{j=1}^{J} c^{(j)} C^{(j)}(x, \dot{x}), \tag{4}$$

where $C^{(j)}(x, \dot{x})$ is a candidate function with $c^{(j)}$ as its scalar coefficient. Substituting in Eq. 4 into Eq. 3 results in

$$\mathbf{Qc} = \mathbf{r}, \tag{5}$$

where

$$\mathbf{c} = \left[c^{(1)}, c^{(2)}, \dots, c^{(J)}\right]^T, \tag{6}$$

$$\mathbf{r} = \left[r^{(1)}, r^{(2)}, \dots, r^{(P)}\right]^T, \tag{7}$$

and the elements of the matrix $\mathbf{Q}$ are computed as

$$\mathbf{Q}^{(p,j)} = \int_{\gamma^{(0)}}^{\gamma^{(p)}} \dot{x} C^{(j)}(x, \dot{x}) dt \,. \tag{8}$$

Equation 5 is a linear system of equations that can be solved to compute the coefficients $c^{(j)}, j = 1,2, \dots, J$. The result is a mathematical model for the nonconservative (damping) force acting on the oscillator.

Phase II identifies a model for the elastic force, $K(x)$, by either using Lagrange's equation applied to a modified Lagrangian [27] computed using the mechanical and kinetic energies or using Newton's second law [28] to solve for $K(x)$ directly. We will use the second approach in the proposed methodology, which simply requires rewriting the equation of motion to be

$$K(x) = F(t) - m\ddot{x} - B(x, \dot{x}), \tag{9}$$

where the model for $B(x, \dot{x})$ determined in Phase I is used. At this point, we propose a model for the conservative force of the form

$$K(x) = \sum_{l=1}^{L} a_l A_l(x) \,, \tag{10}$$

where each $A_l(x)$ is a candidate function with an unknown coefficient $a_l$. The unknown coefficients are determined using curve fitting techniques, linear regression, or constrained optimization. The outcome of EDDI is a mathematical model for the dynamics of the nonlinear SDOF oscillator produced entirely from measured response data.

### 2.2. Wavelet-Bounded Empirical Mode Decomposition (WBEMD)

The EMD introduced by Huang et al. [32] is a signal decomposition technique that extracts a finite basis of nearly orthogonal, monochromatic IMFs from nonlinear, oscillatory signals. Each IMF is associated with a characteristic time scale and are often representative of individual nonlinear normal modes [47], defined to be periodic oscillations that are not necessarily synchronous. EMD relies on a sifting algorithm to decompose an oscillatory signal, $y(t)$, into $N$ nearly orthogonal IMFs, such that

$$y(t) = \sum_{i=1}^{I} y_i(t) + e_{I+1}(t), \qquad \max(|e_{I+1}(t)|) < \tau, \tag{11}$$

where $e_{I+1}(t)$ is a remainder signal that has an amplitude that is less than the tolerance $\tau$. The goal of EMD is to extract IMFs that are physically and mathematically representative of individual harmonics contained in the original signal and that reconstruct the original signal when summed.

Even though IMFs are known to be mathematically and physically significant [47–49], the EMD suffers from a few issues that affect the usefulness and meaning of the IMFs. First, EMD often extracts more IMFs than harmonic components present in the original signal, such that the analyst must carefully select only the IMFs that correspond to physically present components [32,50,51]. A second issue is that of *mode mixing* where multiple harmonic components are mixed into a single IMF, such that the IMF loses physical significance. Mode mixing is caused either due to intermittency or due components being closely spaced in frequency. Intermittency occurs when a component dies out before the end of the signal, which occurs frequently in transient nonlinear vibrations (e.g., free-response). In the context of EMD, closely spaced implies the frequency ratio between two components is $\leq 2$, which also frequently occurs in vibrating structures. While many techniques have been developed to address these two issues (and others), we will only discuss WBEMD because it is directly used in the proposed methodology.

WBEMD [34,35] solves the issues of spurious IMFs and mode mixing by ensuring that each IMF is isolated around and associated with a specific harmonic component contained in the original signal. The analyst first determines the characteristic frequencies, $\omega_i$, of the harmonic components contained in the signal using the maximum wavelet transform [34]. WBEMD uses an optimization routine to systematically extract IMFs that are each localized around a single characteristic frequency. The optimization routine starts by first applying EMD to the signals

$$u_+(t) = z_i(t) + s(t), \tag{12a}$$

$$u_-(t) = z_i(t) - s(t), \tag{12b}$$

where $z_i(t)$ represents the current remainder signal defined as

$$z_i(t) = x(t) - \sum_{q=1}^{i-1} x_i(t)\,, \tag{13}$$

$s(t)$ is a masking signal defined as

$$s(t) = \max[z_i(t)]\ (\alpha \sin(\beta \omega_i t)), \tag{14}$$

and $\alpha$ and $\beta$ are free parameters used by the optimization routine. The use of a masking signal was introduced by Deering and Kaiser [52] and it simultaneously eliminates mode mixing due to intermittency and closely spaced frequencies because it fools EMD into mixing the desired component with the masking signal instead of other components. The masking signal is then removed by averaging the two IMFs together resulting in only the desired harmonic component. However, the difficulty arises in selecting the optimal masking signal to extract the desired component, which is the core problem that WBEMD solves.

After averaging out the masking signal, the resulting IMF is transformed to the maximum wavelet domain [34] where a bounding function is used to measure the isolation of the IMF around the characteristic frequency. The bounding function is constructed using a Gaussian distribution and the area underneath it provides the measure of the isolation of the IMF. A poorly isolated IMF will result in a bounding function that is wide and has a large area, while a well-separated IMF results in a bounding function that is narrow and has a small area. WBEMD leverages this observation to extract well-separated IMFs by optimizing the parameters $\alpha$ and $\beta$ until the area of the bounding function reaches a global minimum. The outcome of this optimization is a well-separated IMF that is representative of a single harmonic component contained in the original signal. WBEMD is available for download at [53]. One key difference of WBEMD compared to EMD is that when the current remainder signal is monochromatic, it is chosen as the final IMF instead of applying WBEMD again. The reason for this is that WBEMD is intended to extract IMFs that are physically representative of harmonic components contained in the signal, such that applying EMD (or WBEMD) to a monochromatic remainder signal results in an IMF and remainder that are not physically meaningful. Thus, the original signal can be reconstructed exactly by summing the IMFs directly when extracted using WBEMD.

### 2.3. The Proposed DEDDI Method

DEDDI synergistically combines WBEMD and EDDI to produce an energy-based, data-driven methodology for nonlinear system identification of MDOF systems undergoing nonlinear oscillations. Specifically, DEDDI leverages WBEMD to extract nearly orthogonal, monochromatic IMFs from the displacement response of each DOF, each representing a single harmonic component contained in the signal. The benefit of this decomposition is that the IMFs can be treated as the response of an SDOF oscillator, such that EDDI can be applied to each IMF sequentially. We note that this treatment is similar to the construction of intrinsic modal oscillators [54,44,46], except that we do not make any assumptions about the linearity of the oscillators that the IMFs represent. Thus, unlike intrinsic modal oscillators, the nonlinear modal interactions captured by the IMFs are not treated like external forcing. A flowchart for the DEDDI process is provided in Fig. 1 and we explain the method in detail in the following.

DEDDI begins with the set of measured displacements for each DOF in the system. While DEDDI can be applied in cases where only some of the DOFs are measured (e.g. continuous structures), the present formulation results in a mathematical model for only the measured DOFs. Thus, if applied to a continuous system, DEDDI will produce a reduced-order model for the dynamics of the system in the form of a set of coupled ordinary differential equations. Applications to continuous systems is left open for future research and this work will focus on discrete MDOF systems. We note that the displacements can be recovered by measuring them directly or by integrating and filtering acceleration or velocity measurements.

The first step of DEDDI is to decompose the displacements into IMFs representative of their harmonic components. For the displacement $x_n(t)$ of the $n$th DOF, WBEMD produces the IMFs $x_{ni}(t), i = 1, \dots, I_n$ where the total number of IMFs may be different for each DOF (e.g., if one DOF is at a node for a particular mode). The original displacement can be recovered by summing the IMFs directly. The second step is to compute the velocities and accelerations of the IMFs, which is accomplished by numerically differentiating them (we use *gradient* in MATLAB® to maintain the size of the vectors). This results in the velocity and acceleration IMFs defined as

$$\dot{x}_{ni}(t) = \frac{dx_{ni}}{dt}, \ddot{x}_{ni}(t) = \frac{d\dot{x}_{ni}}{dt} = \frac{d^2 x_{ni}}{dt^2}. \tag{15}$$

Assuming that the velocity and acceleration time series are available, one could apply WBEMD to those signals instead to compute the velocity and acceleration IMFs. However, we found that using the numerical derivatives of the displacement IMFs produced better results than computing the IMFs of the velocity and acceleration time series directly. To avoid confusion, we will use displacement IMF, velocity IMF, or acceleration IMF when we refer to a specific signal. We will use IMF to refer to the component in general, typically in the context of treating it as an SDOF oscillator.

With the IMFs and their derivatives computed, DEDDI proceeds by treating each IMF as an SDOF oscillator and applies EDDI to each of them directly. The third step is to compute the kinetic energy of each IMF using the mass of the $n$th DOF, such that

$$T_{ni}(t) = \frac{1}{2} m_n \dot{x}_{ni}^2(t). \tag{16}$$

We note that the nearly orthogonal nature of the IMFs allows us to compute these individual kinetic energies and neglect the mixed terms that arise when computing the kinetic energy from an expanded velocity. Specifically, if one expands the velocity as

$$\dot{x} = \sum_{i=1}^{I} \dot{x}_i \, , \tag{17}$$

then the kinetic energy for that DOF is

$$T(t) = \frac{1}{2} m \dot{x}^2 = \frac{1}{2} m \left( \sum_{i=1}^{I_n} \dot{x}_i \right)^2 = \frac{1}{2} m \sum_{i=1}^{I} \dot{x}_i^2 + m \sum_{i=1}^{I} \sum_{j=1}^{I} \dot{x}_j \dot{x}_i \, . \tag{18}$$

Due to the nearly orthogonal nature of the IMFs, the mixed term at the end of the Eq. (**Error! Bookmark not defined.**) will be nearly zero and small enough to neglect. For discrete systems, the mass of each DOF can be directly measured, but the mass that should be used for continuous systems is less obvious. One choice is to use a lumped mass approximation by dividing the structure into segments centered around each, but this remains an open question for future research.

The fourth step is to determine the time instants, $\gamma_{ni}^{(p)}, p = 0,1,\ldots,P_{ni}$, where each displacement IMF is zero and assume that its corresponding potential energy is also zero at those times. Next, we shift time to begin at $\gamma_{ni}^{(0)}$ and estimate the dissipated energies for the IMF as

$$r_{ni}^{(p)} = T_{ni}\left(\gamma_{ni}^{(0)}\right) - T_{ni}\left(\gamma_{ni}^{(p)}\right). \tag{19}$$

The estimated dissipated energies can be plotted against time if desired, which can provide insight into the type of dissipation mechanisms that arise in the dynamics of the IMF.

The sixth step is to propose a mathematical model that is linear in parameters for the dissipated energy of each IMF of the form

$$B_{ni}(x_{ni}, \dot{x}_{ni}) = \sum_{j=1}^{J} c_{ni}^{(j)} C_{ni}^{(j)}(x_{ni}, \dot{x}_{ni}) \, , \tag{20}$$

where $c_{ni}^{(j)}$ is an unknown coefficient and $C_{ni}^{(j)}(x_{ni}, \dot{x}_{ni})$ represents a single function used for the model. Next, we compute the model dissipated energy for each model term and construct the linear system of equations as

$$\mathbf{Q}_{ni}\mathbf{c}_{ni} = \mathbf{r}_{ni}, \tag{21}$$

where the elements of the matrix $\mathbf{Q}_{ni}$ are computed as

$$Q_{ni}^{(p,j)} = \int_{\gamma_{ni}^{(0)}}^{\gamma_{ni}^{(p)}} \dot{x}_{ni} C_{ni}^{(j)} \, (x_{ni}, \dot{x}_{ni}) dt \,, \tag{22}$$

and

$$\mathbf{c}_{ni} = \left[c_{ni}^{(1)}, c_{ni}^{(2)}, \dots, c_{ni}^{(J)}\right], \tag{23}$$

$$\mathbf{r}_{ni} = \left[r_{ni}^{(1)}, r_{ni}^{(2)}, \dots, r_{ni}^{(P)}\right]. \tag{24}$$

The coefficients $\mathbf{c}_{ni}$ can be computed using linear regression (i.e., $\mathbf{c}_{ni} = \mathbf{Q}_{ni}^{-1}\mathbf{r}_{ni}$ solved in a least-squares sense), curve fitting, or other optimization methods as needed. In this work, we use least-squares regression to compute the coefficients.

With the damping model coefficients computed, next we compute the time-series for the dissipative force of the $i$th IMF of the $n$th DOF as

$$B_{ni}(t) = \sum_{j=1}^{J} c_{ni}^{(j)} C_{ni}^{(j)}(t) \,. \tag{25}$$

We note that at this point we are treating $B_{ni}(t)$ as a time series and no longer as an unknown function of the displacement and velocity IMFs. Next, we compute the conservative force for the IMF by leveraging Newton's second law, such that

$$K_{ni}(t) = -m_n \ddot{x}_{ni}(t) - B_{ni}(t), \tag{26}$$

where the mass of the DOF is used to compute the inertial force.

The application of EDDI concludes at this point with the output being the dissipative and conservative forces acting on each IMF. Next, we leverage the fact that the IMFs can be summed to produce the original time series to obtain estimates for the internal forces acting on each DOF. Specifically, we sum the nonconservative and conservative forces from each IMF to produce the corresponding force at each DOF:

$$B_n(t) = \sum_{i=1}^{I_n} B_{ni}(t) \,, K_n(t) = \sum_{i=1}^{I_n} K_{ni}(t) \,. \tag{27}$$

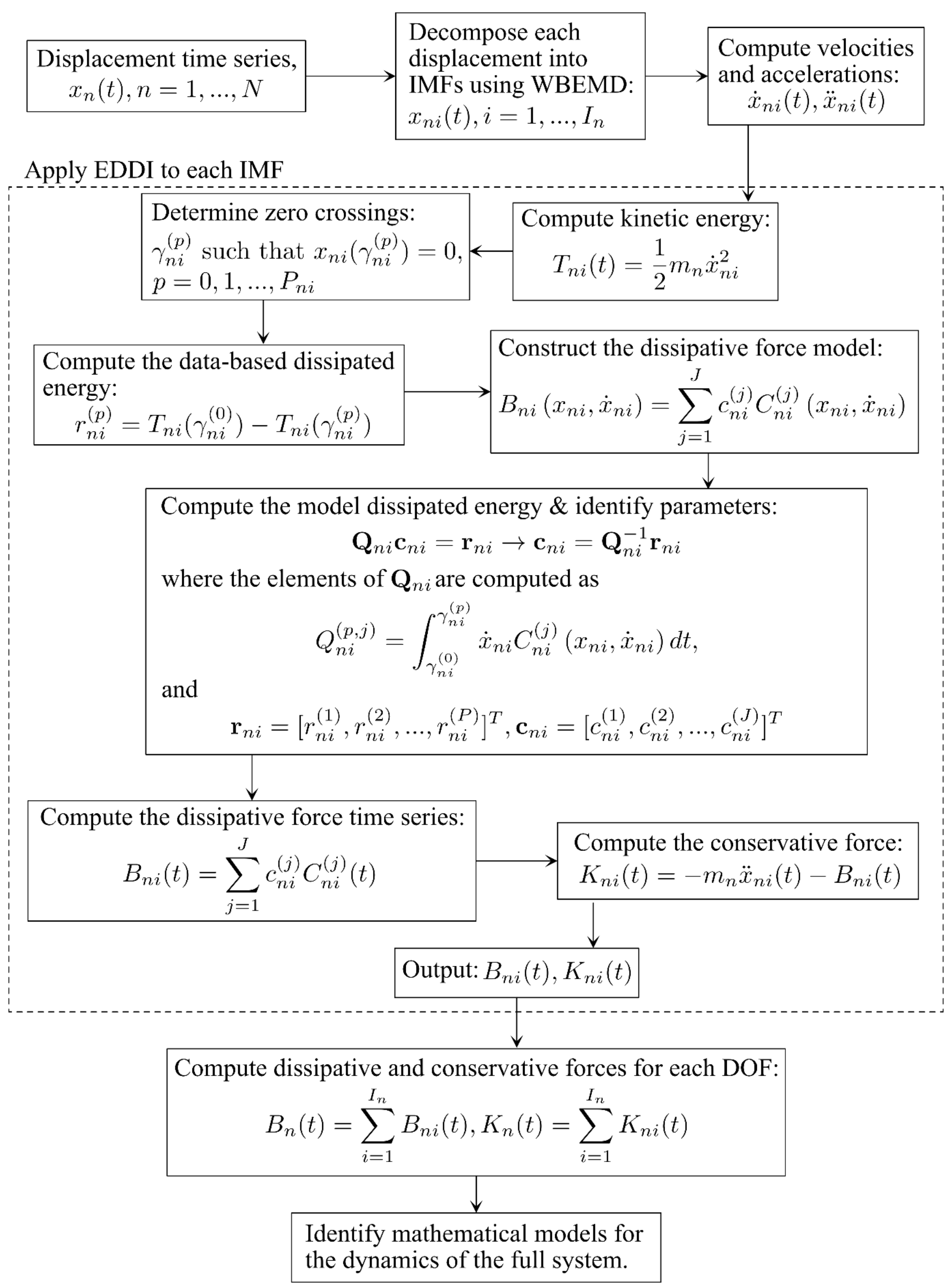

Figure 1. Flowchart depicting the process for DEDDI with the application of EDDI to each IMF detailed.

The result of this is the internal nonconservative and conservative forces acting on each DOF, which can be used to identify a model for the dynamics of the system. The advantage of this process is that it separates the nonconservative and conservative forces, such that one can identify models for each of them separately. This allows one to vary the model for the damping and stiffness separately, while ensuring that they don't influence each other. At this stage, we propose mathematical models for the internal forces, then identify the unknown coefficients. The identification can be performed using a variety of methods including sparse regression or any curve-fitting procedure. In this work, we will leverage constrained optimization to maximize the R-squared value computed between the estimated forces and the model forces. A detailed discussion of this implementation is provided in the next section when it is applied to the example problem.

## 3. Experimental Demonstration of DEDDI

This section focuses on the application and demonstration of DEDDI to a two-story (two DOF) tower with nonlinear stiffness connecting the two floors. We start with a description of the experimental system, the measurement setup, and the measurement results. We then apply DEDDI, walk through every step in sequence, and provide the identified parameters produced by DEDDI, SIVA, and SINDy. Lastly, we validate the identified model by simulating the response of the structure for the measurement case used in the identification and compare with the results from SIVA and SINDy. We also validate the model produced using DEDDI using measurement cases not used in the identification.

### 3.1. Experimental System, Measurement Setup, and Measurement Results

The experimental setup consists of a two-story tower (Fig. 2) composed of a linear oscillator (LO) and a nonlinear oscillator (NO) that was previously studied in [55,56,13]. The LO is installed on an optical table using a pair of steel flexures that act as linear springs, while the NO is coupled to the LO through a pair of steel flexures and a pair of thin steel wires. The coupling flexures have a large rectangular hole cut in them, introduce a linear coupling stiffness, and restrict the NO to horizontal motion. The steel wires are clamped to the LO on their ends and to the NO in their centers using 10-32 set screws, which causes them to introduce a strong geometric nonlinearity that couples the two masses. The LO and NO are both made from aluminum, and each have a mass of 1.37 kg. Note that additional masses were installed on the NO to ensure its mass was equal to that of the LO. Readers are advised to check [55,56] for full details, materials, and dimensions of all components.

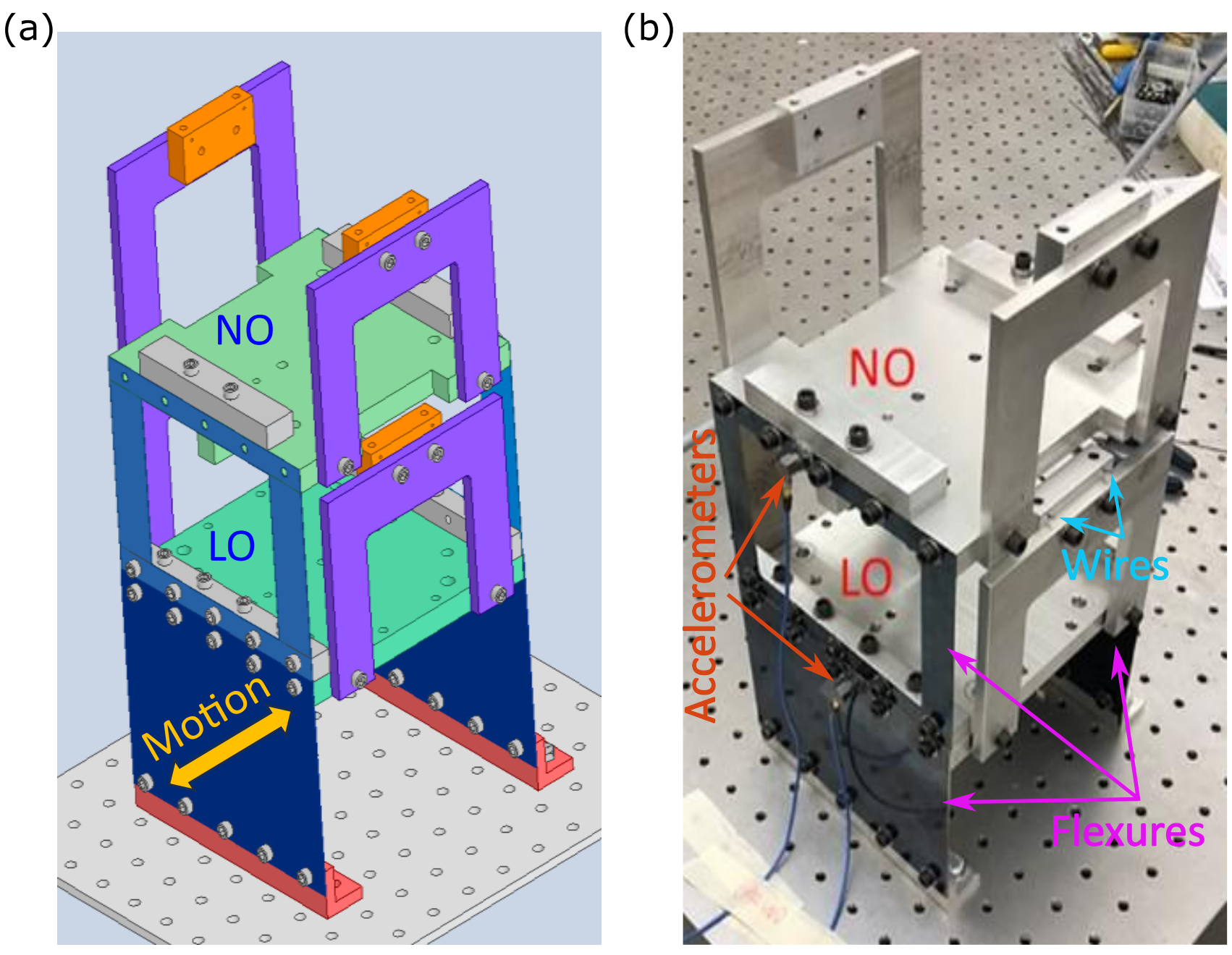

Figure 2. (a) CAD model depicting the experimental system composed of two oscillators. (b) The experimental system with two floors of equal mass with the instrumentation shown.

The experimental measurements used in this work were first reported in [55] and we summarize them here. The system was excited by hitting the NO using a PCB Piezotronics modal hammer (model 086C03). The resulting motion was measured using two PCB accelerometers (model 353B15) with nominal sensitivities of 1 mV/(m/s$^2$). An Abacus 906 data acquisition system (Data Physics, San Jose, CA USA) and the Data Physics software suite was used to record the data with a sampling rate of 16,384 Hz for 8 seconds. The velocities were obtained by numerically integrating the accelerations and high-pass filtering the output using a third-order Butterworth filter with a cut-off frequency of 2 Hz. This process was applied to the velocities to produce the displacements.

We present the measured impact force and displacements for the measurement case that will be used in DEDDI in Fig. 3a and b, respectively. The measured force shows that a clean input without any double impacts. The displacement time series shows several beats in the response of the LO and one pronounced beat for the NO. These beats arise due to the strongly nonlinear interactions that arise between the LO and NO, typically through a 3:1 internal resonance [55] with the measured displacements, we include the continuous wavelet transforms (CWT) for both the LO and NO in Fig. 3b, which provide a time-frequency representation of the signal content. The CWT spectra are normalized to be on the same scale (from 0 to 1) and shaded such that darker shading represents higher energy content at a given time and frequency. The CWT spectra reveal that two harmonic components are clearly visible for the LO, but the NO is dominated by a low-frequency component. The second component is present in the displacement of the NO, but it is substantially smaller than the first component and is not visible in the CWT spectrum. The two components observed in the CWT spectra have frequencies of around 8 and 20 Hz, and these frequencies will be used for the characteristic frequencies in the decomposition using WBEMD in the next subsection.

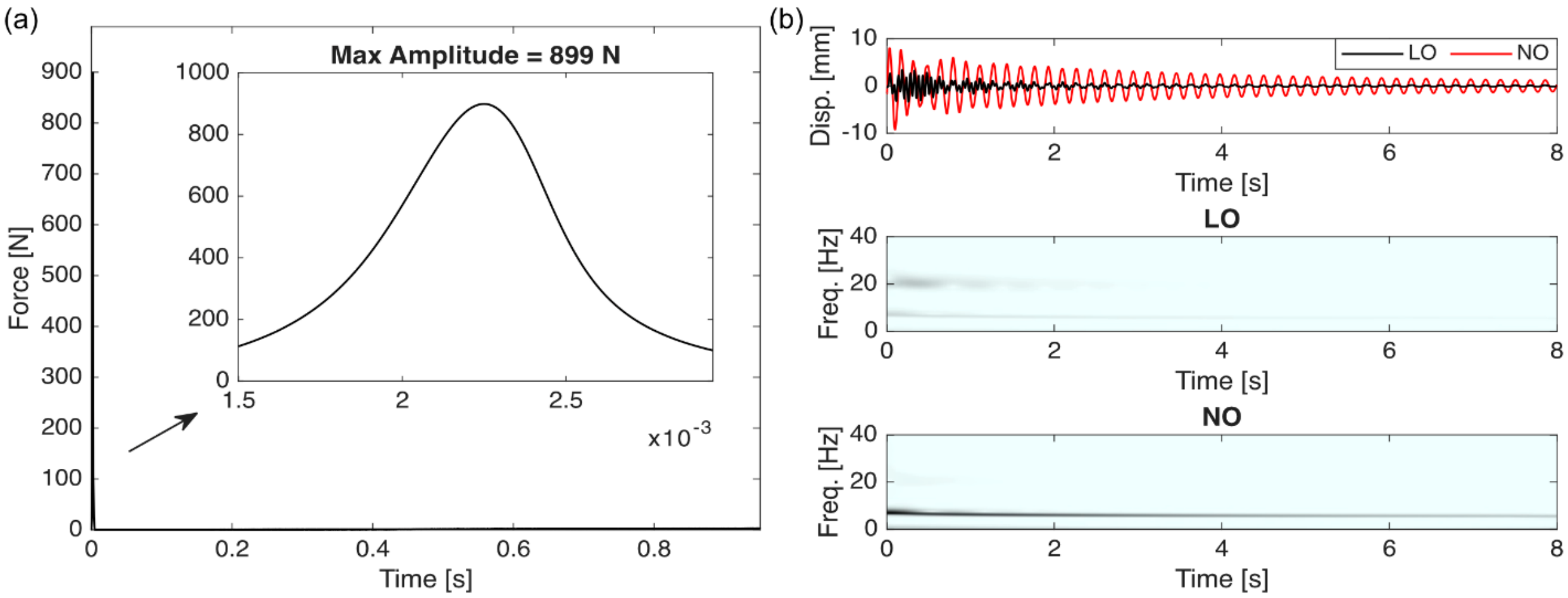

Figure 3. The measured displacement response of the system for an impact of 899 N with the time series, continuous wavelet transforms, and the maximum wavelet transform shown for each DOF.

### 3.2. Application of DEDDI

Using the measured response shown in Fig. 3, we start by applying WBEMD to decompose each displacement signal into two IMFs with the characteristic frequencies of 8 and 20 Hz. The resulting IMFs are shown in Fig. 4a and b for the LO and NO, respectively. Here, we label the LO as the first DOF indicated by a 1 in the leading subscript and the NO as the second DOF indicated by a leading subscript of 2. The trailing subscripts correspond to the IMF number ordered from lowest to highest frequency to match the two nonlinear normal modes in the system. We present the time series, CWT spectrum, and the instantaneous frequency for each IMF in Fig. 4. The instantaneous frequency is shown as the solid line in the CWT spectra and is computed using the Hilbert transform. The time-series reveal that the amplitudes of the IMFs are slowly varying and do not monotonically decay as would be seen in a linear response. For example, the amplitude of IMF $x_{12}(t)$ (the second IMF of the LO displacement) exhibits many nonlinear beats that is indicative of an internal resonance in the response [46] (in this case, this arises through a 3:1 internal resonance between the first and second nonlinear normal modes [55]). The CWT spectra and instantaneous frequencies confirm that the IMFs are well separated, and each represents a single component contained in the original signals. Most importantly, each IMF is monochromatic, such that we can treat each of them as a SDOF oscillator and apply EDDI to estimate their internal nonconservative and conservative forces.

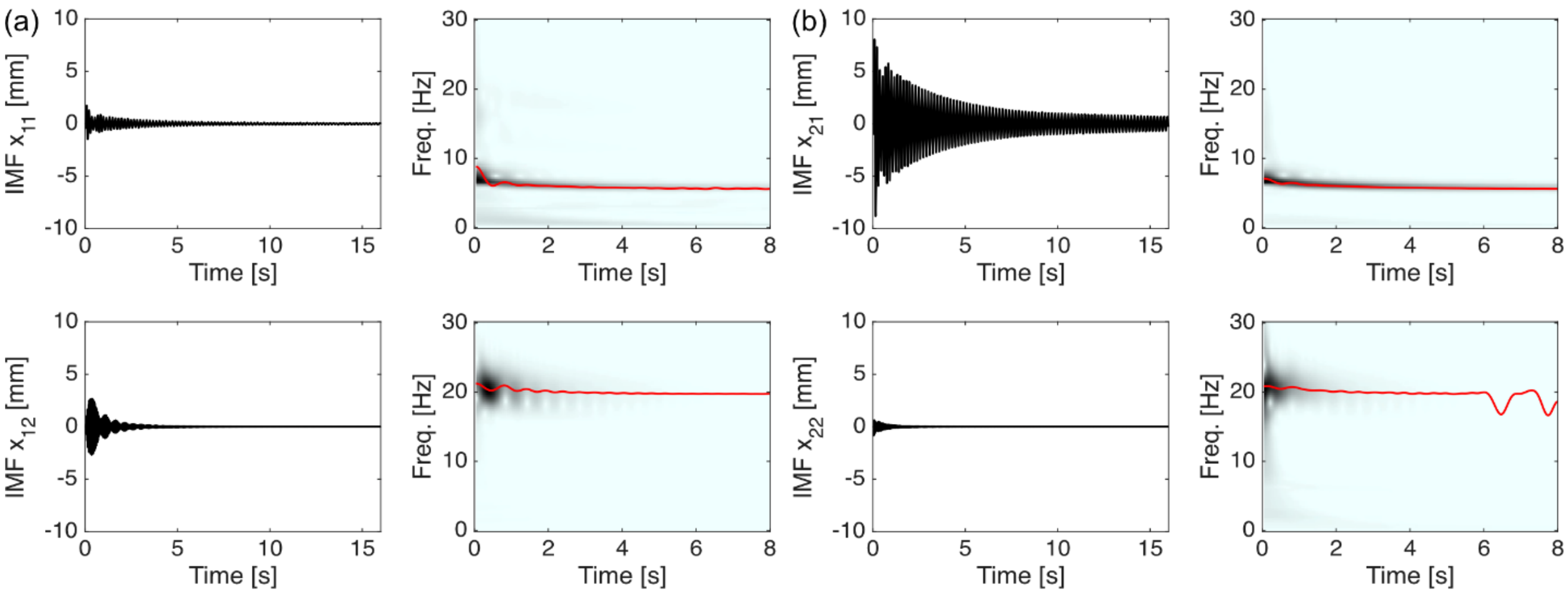

Figure 4. The time series and CWT for the IMFs extracted from the displacement of (a) the LO and (b) the NO. The instantaneous frequency is plotted on top of the CWT to show that the signals are monochromatic.

We proceed by computing the IMF velocities and accelerations numerically using the *gradient* function in MATLAB, which preserves the size of the vectors. Using the IMF velocities, we compute their corresponding kinetic energies and present them in the first row of Fig. 5. Note that the kinetic energies are plotted on different scales and are not quantitatively comparable as presented here. Interestingly, the kinetic energies also exhibit slowly varying amplitudes and do not follow a monotonic decay as is seen for linear systems. The slow changes in the amplitudes of the kinetic energies are related to the sharing of energy between the nonlinear normal modes in the system and may provide a path to determining how much energy is shared between them. This remains an open question for future research and investigation. Next, we determine the times where the IMF displacements are zero and assume that the potential energy of each IMF is also zero at those times. These times are computed by first finding the time instants where the displacement crosses zero, then linear interpolating to find the time where the displacement is zero. We plot the kinetic energy at those times as circles on the first row of Fig. 5.

Next, we shift time to begin at the first instant where the IMF displacement is zero, then compute the data-based dissipated energy using Eq. 19. This computation is performed separately for each IMF as the times of zero IMF displacement are different for each IMF. We present the dissipated energy for each IMF in the middle row of Fig. 5 as black circles. In general, the dissipated energies increase rapidly before asymptotically approaching a maximum value; however, they do not increase monotonically as would be expected in a linear system. Instead, oscillations arise in the dissipated energies of all IMFs, though these are most pronounced in the first three. The oscillations observed in the dissipated energies of the IMFs are likely related to the energy transfer that arises between the components through the nonlinearity. More importantly, the oscillations are secondary to the overall increase in the dissipated energies, such that they do track the overall loss of energy over time.

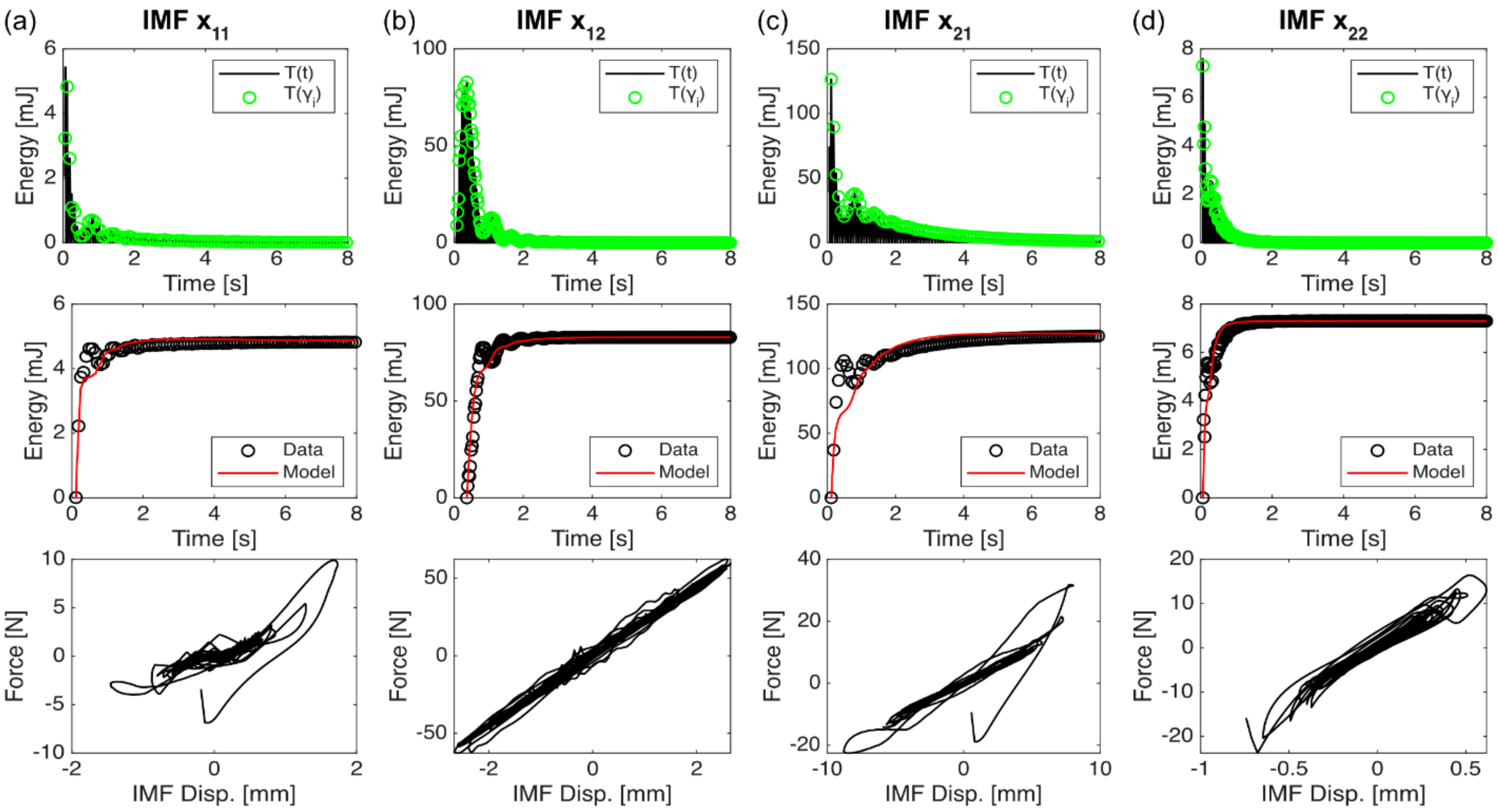

Figure 5. The kinetic energies, dissipated energies, and conservative force for each IMF.

At this point, we propose a model for the nonconservative force for each IMF and use those to estimate the unscaled dissipated energies, leaving the coefficients unknown for identification. In this case, we select the following model for the nonconservative forces:

$$B_{ni}(x_{ni}, \dot{x}_{ni}) = c_{ni}^{(1)}\dot{x}_{ni} + c_{ni}^{(2)}x_{ni}^2\dot{x}_{ni}, \quad n, i = 1,2, \tag{28}$$

where the first term represents a linear viscous damper and the second term arises due to a viscous damper moving transversely [57]. We note that one may be able to capture the energy transfers by incorporating model terms that couple the IMFs, though this is left open for future research at this time. Using the model, we construct the unscaled dissipated energies, then form the linear system of equations defined in Eq. 21. We compute the unknown coefficients $c_{ni}^{(p)}$ using least squares regression and provide their values in Table 1. Using the identified coefficients, we compute the dissipated energies and plot them as the solid lines on the middle row of Fig. 5. These comparisons show that the identified models reproduce the overall behavior of the dissipated energies, with the model for IMF $x_{21}$differing the most from the data.

Table 1. The values of the damping parameters identified for each IMF.

| **Parameter** | **IMF $x_{11}$** | **IMF $x_{12}$** | **IMF $x_{21}$** | **IMF $x_{22}$** |
|---|---|---|---|---|
| $c_{ni}^{(1)}$ [Ns/m] | $-0.8681$ | 3.366 | $-0.2397$ | 0.4656 |
| $c_{ni}^{(2)}$ [Ns/m$^3$] | $2.797 \times 10^7$ | $1.340 \times 10^6$ | $3.705 \times 10^5$ | $1.530 \times 10^8$ |

At this point, we compute the time series for the nonconservative force for each IMF, then leverage those to compute the conservative force time series using Eq. 26 and the IMF accelerations. The resulting conservative forces are depicted in the last row of Fig. 5 as restoring forces by plotting them as functions of the IMF displacement. Interestingly, the conservative forces for IMFs $x_{11}$ and $x_{21}$ exhibit the most nonlinear behavior while those for IMFs $x_{12}$ and $x_{22}$ are primarily linear. This is likely a reflection of the strong linear grounding stiffness and the fact that the second nonlinear normal mode is dominated by motion in the LO [55]. At this point, the application of EDDI to each IMF concludes as no model is needed for the conservative force. Instead, we sum the nonconservative and conservative forces for the two IMFs from each DOF to produce the internal forces acting on each DOF. We depict the internal nonconservative and conservative forces acting on the LO and NO as the solid lines in Figs. 6a and b, respectively. Thus, the combination of WBEMD and EDDI allows one to convert measured transient responses into the internal forces acting on each DOF. Crucially, this approach isolates the nonconservative from the conservative forces, so that one can identify models for them separately and without contamination.

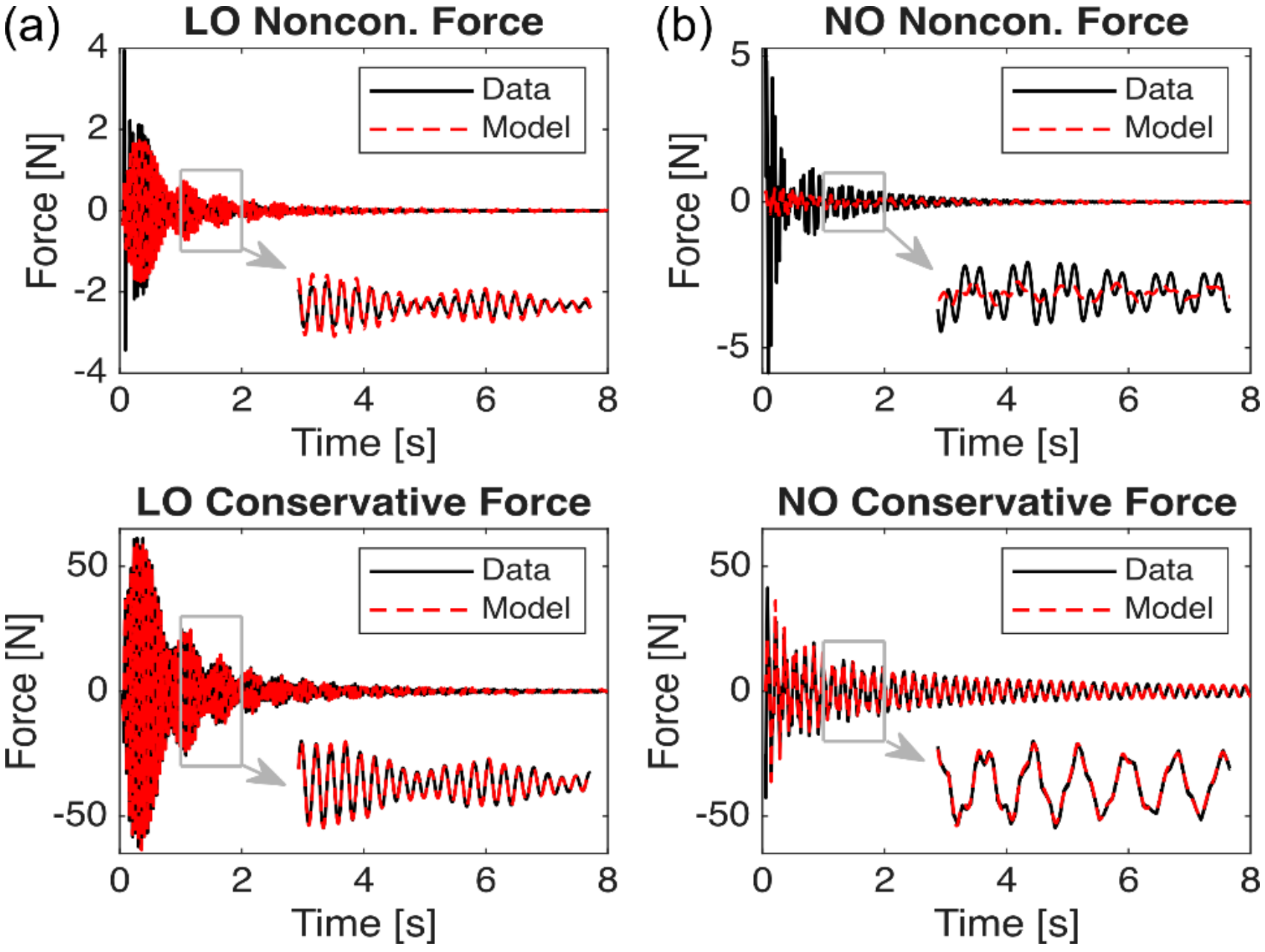

Figure 6. The total nonconservative and conservative force for (a) the LO and (b) the NO.

We proceed by selecting models for the damping and stiffness of the system. While one can choose any library of functions for the models, some understanding of the physical structure can help eliminate terms in that library. For example, the fact that the NO is not coupled to the ground directly in the experiment allows one to select only coupling terms for the models of the nonconservative and conservative forces of the NO. In this application, we select the same model structure that was used in the original papers on the LO-NO system [55,56] as well as in the original paper on SIVA [13]. The models for the internal nonconservative forces are

$$B_1(\dot{x}_1, \dot{x}_2) = b_g \dot{x}_1 + b_c(\dot{x}_1 - \dot{x}_2), \tag{29a}$$

$$B_2(\dot{x}_1, \dot{x}_2) = b_c(\dot{x}_2 - \dot{x}_1), \tag{29b}$$

$$K_1(x_1, x_2) = k_g x_1 + k_1(x_1 - x_2) + k_3(x_1 - x_2)^3, \tag{29c}$$

$$K_2(x_1, x_2) = k_1(x_2 - x_1) + k_3(x_2 - x_1)^3, \tag{29d}$$

where the subscripts "1" corresponds to the LO and "2" to the NO. Note that we have chosen to use a subscript of "3" for the nonlinear coupling stiffness due to the use of a cubic nonlinearity, such that $k_2$ would correspond to the stiffness of a quadratic term if one were included. We also note that the models for the nonconservative forces used here are different than those used in modeling the dissipation of the IMFs. We originally incorporated a nonlinear damping term of the form $b_2(\dot{x}_2 - \dot{x}_1)(x_2 - x_1)^2$, but found that the resulting parameter value was relatively weak and the reproduction of the measurements was better when that term was excluded. Furthermore, we chose to exclude that term to match the models used in [55,56,13] to form a fair comparison.

The unknown parameters were identified using constrained global optimization to maximize the R-squared value between the forces estimated using DEDDI and those produced by the model. We identified the damping parameters separately from the stiffness parameters, but simultaneously for each oscillator. The objective function for the damping identification maximized the combined R-squared values for the nonconservative forces of LO and NO, while that for the stiffness coefficients combined the R-squared values for the conservative forces. To perform the optimization, we used *patternsearch* in MATLAB® with the maximum function evaluations set to $10^9$, the max iterations to $10^9$, the mesh tolerance to $10^{-12}$, the step tolerance to $10^{-14}$, and all other options set to their default values. These values were chosen to force *patternsearch* to find a true global maximum and to ensure that it only stopped when the mesh tolerance dropped below $10^{-12}$. Similar settings were used in the time-series optimization method applied in [55]. The initial guesses, lower bounds, and upper bounds for each parameter are provided in Table 2, and we selected based on experience and to provide a wide range for the optimization to explore. The identified values produced by DEDDI are also provided in Table 2 along with the values identified by SIVA and time-series optimization for comparison. Lastly, we report the combined mean squared error (MSE) computed between the measured and simulated displacements for each oscillator using the given model, which shows that DEDDI produces the lowest MSE value. Further validation of DEDDI is provided in the next subsection.

Using the identified parameters, we compute the time series for the model nonconservative and conservative forces, then plot them in Figs. 6a and b for the LO and NO, respectively. There is strong agreement between the conservative forces with the R-squared values being 0.9891 and 0.9214 for the LO and NO, respectively. However, the agreement in the nonconservative forces is weaker for the LO with an R-squared value of 0.6323 and significantly weaker for the NO with an R-squared value of 0.1125. For the LO, we find that the model force captures the overall profile, including the beats, but its amplitude is higher than that of the experiment. The model also does worse at the very beginning, though the features in the data may be caused by end effects introduced by WBEMD. In contrast, the model nonconservative force for the NO does not accurately reproduce the force estimated using DEDDI in either amplitude or overall shape. This disagreement suggests that a better model for the damping of NO may exist, and if desired, one could leverage this data to explore alternative models. In the present work, we proceed with the chosen model to allow for fair comparison with other methods as discussed previously.

Table 2. The initial guesses, lower bounds, and upper bounds used in DEDDI and the identified values for the damping and stiffness coefficients using DEDDI, SIVA [13], and the time-series optimization used in [55].

| Parameter | Initial Guess | Lower Bound | Upper Bound | DEDDI | SIVA | Time-Series Optimization [55] |
|---|---|---|---|---|---|---|
| $b_g$ [Ns/m] | 2 | 0 | 10 | 4.2142 | 3.8405 | 2.177 |
| $b_c$ [Ns/m] | 2 | 0 | 10 | 0.7458 | 0.7856 | 1.075 |
| $k_g$ [N/m] | $10^4$ | $-10^5$ | $10^5$ | 18799 | 19273 | 19239 |
| $k_1$ [N/m] | 2000 | $-10^4$ | $10^4$ | 1943.8 | 1947.6 | 1980.1 |
| $k_2$ [N/m$^3$] | $5 \times 10^7$ | $-10^9$ | $10^9$ | $2.985 \times 10^7$ | $2.764 \times 10^7$ | $2.539 \times 10^7$ |
| MSE [mm$^2$] | – | – | – | 0.9132 | 6.616 | 6.564 |

Lastly, before moving to validation of the models, we applied SINDy to identify the equations of motion for the LO-NO system. We used $[\dot{x}_1, \dot{x}_2, x_1, x_2, (x_2 - x_1)^3]$ for the dictionary of model terms, then applied the standard SINDy algorithm to identify the unknown parameters. We originally included the quadratic term $(x_2 - x_1)^2$, but found that this term was reduced to zero and caused the NO damping to become negative. Thus, we chose to remove this term from the library instead of retaining it, which gives SINDy the best opportunity to reproduce the dynamics and makes the resulting model comparable with the others. To apply sparsity conditions, we scaled each parameter by the upper bound used in the DEDDI optimization (Table 2), then removed them if the scaled value was below a threshold of 0.01. all other terms remained. The resulting equations of motion produced by SINDy are provided below:

$$m_1\ddot{x}_1 + 3.1408\dot{x}_1 - 0.6180\dot{x}_2 + 20859x_1 - 1796.9x_2 + 3.649 \times 10^7(x_1 - x_2)^3 = 0, \quad (30a)$$

$$m_2\ddot{x}_2 + 0.2973\dot{x}_1 + 0.683\dot{x}_2 - 1765.9x_1 + 2353x_2 + 3.502 \times 10^7(x_2 - x_1)^3 = F(t). \quad (30b)$$

We find that the equations of motion identified by SINDy lack the symmetry in coupling parameters that arises due to Newton's third law.

### 3.3. Validation of Identified Models

To validate the identified models, we form the governing equation for the two oscillators as

$$m_1\ddot{x}_1 + b_g\dot{x}_1 + b_c(\dot{x}_1 - \dot{x}_2) + k_g x_1 + k_1(x_1 - x_2) + k_3(x_1 - x_2)^3 = 0, \quad (31a)$$

$$m_2\ddot{x}_2 + b_c(\dot{x}_2 - \dot{x}_1) + k_1(x_2 - x_1) + k_3(x_2 - x_1)^3 = F(t), \quad (31b)$$

where all initial conditions are set to zero and $F(t)$ is the experimentally measured impulse. We simulated the response of the system using the experimentally measured force for the identification case using the models identified by DEDDI, SINDy, and SIVA. We present the results for DEDDI in Fig. 7a and b for the LO and NO, respectively. We include comparisons of the displacement time series, CWT spectra, and the Fourier spectra for each DOF. The time series show that the model does a good job of reproducing the measured response for both the LO and NO in both amplitude and phase. The MSE and R-squared values for the time series are provided in each subplot and indicate good agreement. The model displacements agree well with the measured response at both early and late times in the measurements, indicating that both the damping and nonlinearity are well represented. There is a slight shift in the predicted LO displacement, and the NO displacement is a bit lower in amplitude after about 6 seconds. This suggests that the damping in the model may be slightly too high compared to the true dissipation in the experiment. Additionally, the model response also reproduces the beating found in the displacements of both DOFs. The CWT spectra show that the model reproduces the overall signal content, though fewer beats are found in the spectra of the model response for the LO. The two harmonic components are clearly present in the LO response, while the NO is dominated by a single harmonic, as observed in the experiments. Lastly, the Fourier spectra show that the model reproduces the overall signal content.

We present the displacement time series, CWT spectra, and Fourier spectra for simulated results produced by the model identified by SINDy in Figs. 8a and b for the LO and NO, respectively. The time series show that the model produced by SINDy captures the beats in the response of the LO but does not match the overall profile of the oscillations. Additionally, there is a slight phase shift towards the end of the signal, indicating a mismatch in frequency content. The predicted response of the NO does not agree well with the measured response. The amplitude is larger than the experiments, and the beating pattern is also different. There is also a phase shift later in the response, that indicates a frequency mismatch. Looking at the CWT spectra, we find that the SINDY appears to capture the overall phenomena—two harmonics and beating for the LO and only one for the NO with no beating. However, the frequencies of the components are higher than those produced by the experiment, which agrees with the observation of a frequency mismatch from the time series. The comparison of the Fourier spectra indeed shows that the frequency content produced by the SINDy model are shifted to the right. This indicates that the stiffness in the model is too large, and comparing with the other models, it is likely that all stiffness terms are too big.

Next, we consider the results produced by the model identified by SIVA, which are presented in Fig. 8c and d for the LO and NO, respectively. In general, the model agrees with the measurements at early times, then deviates at later times with a phase shift in both the LO and NO displacements. The phase shifts are likely caused by either frequency or damping mismatches in the models. However, the model does a good job of reproducing the beats that appear in both displacements. The predicted amplitude of the LO matches the experiment, but the response of the NO is larger in amplitude at early times, decays to match the experiment around 6 seconds, and finally is lower after 6 seconds. This suggests that the damping in the model is too high, while the stiffness may not be strong enough. Comparing the CWT spectra, we find that the model reproduces the overall signal content and does a better job than DEDDI at reproducing the beats in the second component of the LO. The Fourier spectra reveal that the model reproduces the signal content, but is lower in amplitude at higher frequencies, which further suggests that the effective damping is too large in the model. Overall, the SIVA model reproduces the experiments and is slightly worse than that identified by DEDDI, such that either could be used in practice.

To further validate the model produced by DEDDI, we apply it to simulate the response of the system for lower and higher impact forces. We present the simulated results for an impact of 530 N in Figs. 9a and b for the LO and NO, respectively. Overall, the DEDDI model agrees well with the experiments in the time series, the CWT spectra, and the Fourier spectra. There is a slight amplitude disagreement in the time series, but no phase shift is observed in either displacement. Since the applied force here is relatively close the to force used for identification, we expect the model to do reasonably well. We also simulated the response using the models produced by SIVA and SINDy, but considering brevity, we do not include figures for their results. Instead, we rely on the R-squared values to compare them with DEDDI and provide those in Table 3. Comparing the R-squared values, we find that SIVA performs the best for this force with strong values for both the LO and NO, while SINDy performs the worst.

We now consider the performance of the models for the much higher force of 1986 N, which is more than double the force used in the identification case. We present the displacement time series, CWT spectra, and Fourier spectra for DEDDI in Figs. 10a and b for the LO and NO, respectively. Overall, the predicted response for the LO does not match the experiment at early times but does match at later times. The predicted response of the NO reasonably reproduces the measured response with some amplitude differences and a mismatch in the beating at early times. These results suggest that the nonlinearity in the model is not strong enough and that higher order terms (e.g., $(x_2 - x_1)^5$) should be included if the model is expected to perform well for very high inputs. The CWT spectra show that the content of the NO is well reproduced, while the second component of the LO is not accurately reproduced. Similar mismatches are observed in the Fourier spectra. Comparing the R-squared values for the three models, we find that SIVA again does the best for this case with SINDy being the worst and DEDDI in between the two. These results show that SIVA produces models that perform better over a range of cases than that produced by DEDDI, but DEDDI can produce a better model for a specific measurement case. This outcome is expected because SIVA directly incorporates model validation against additional datasets directly into the identification process, whereas validation must be done separately for DEDDI after the identification.

Table 3. The R-squared values computed between the predicted and measured displacement responses for all three measurement cases for the three models compared.

| **Model** | **LO 530 N** | **NO 530 N** | **LO 899 N** | **NO 899 N** | **LO 1986 N** | **NO 1986 N** |
|---|---|---|---|---|---|---|
| DEDDI | 0.7475 | 0.6833 | 0.7499 | 0.8427 | 0.2974 | 0.5207 |
| SINDy | −0.7108 | −1.009 | −0.8686 | −1.382 | −0.2262 | −0.6422 |
| SIVA | 0.8931 | 0.8765 | 0.5371 | −0.2220 | 0.4816 | 0.7216 |

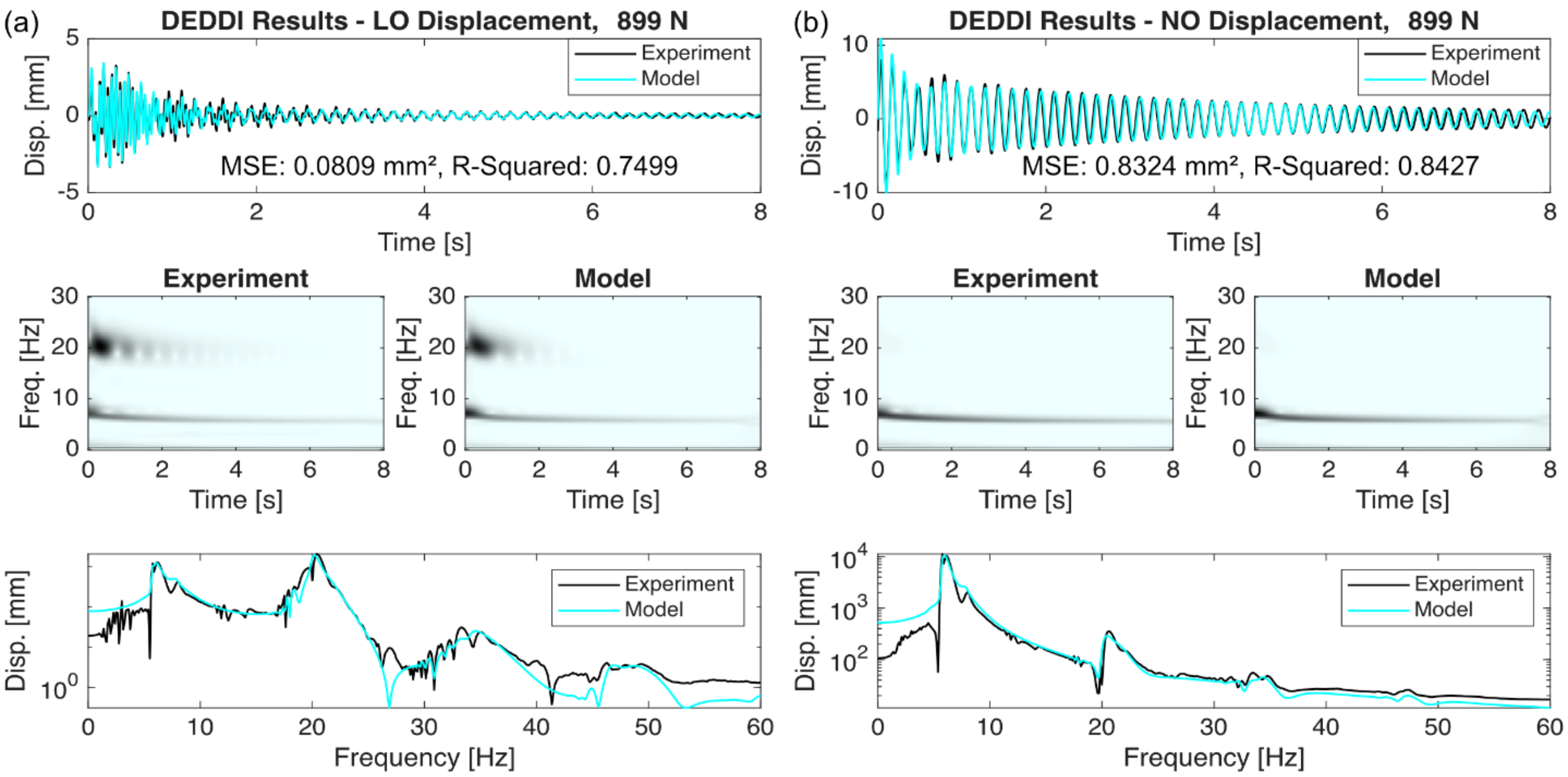

Figure 7. Comparison of the measured and predicted responses using the model produced by DEDDI for the measurement case used for the identification. Time-series, CWT spectra, and Fourier spectra shown for (a) the LO and (b) the NO.

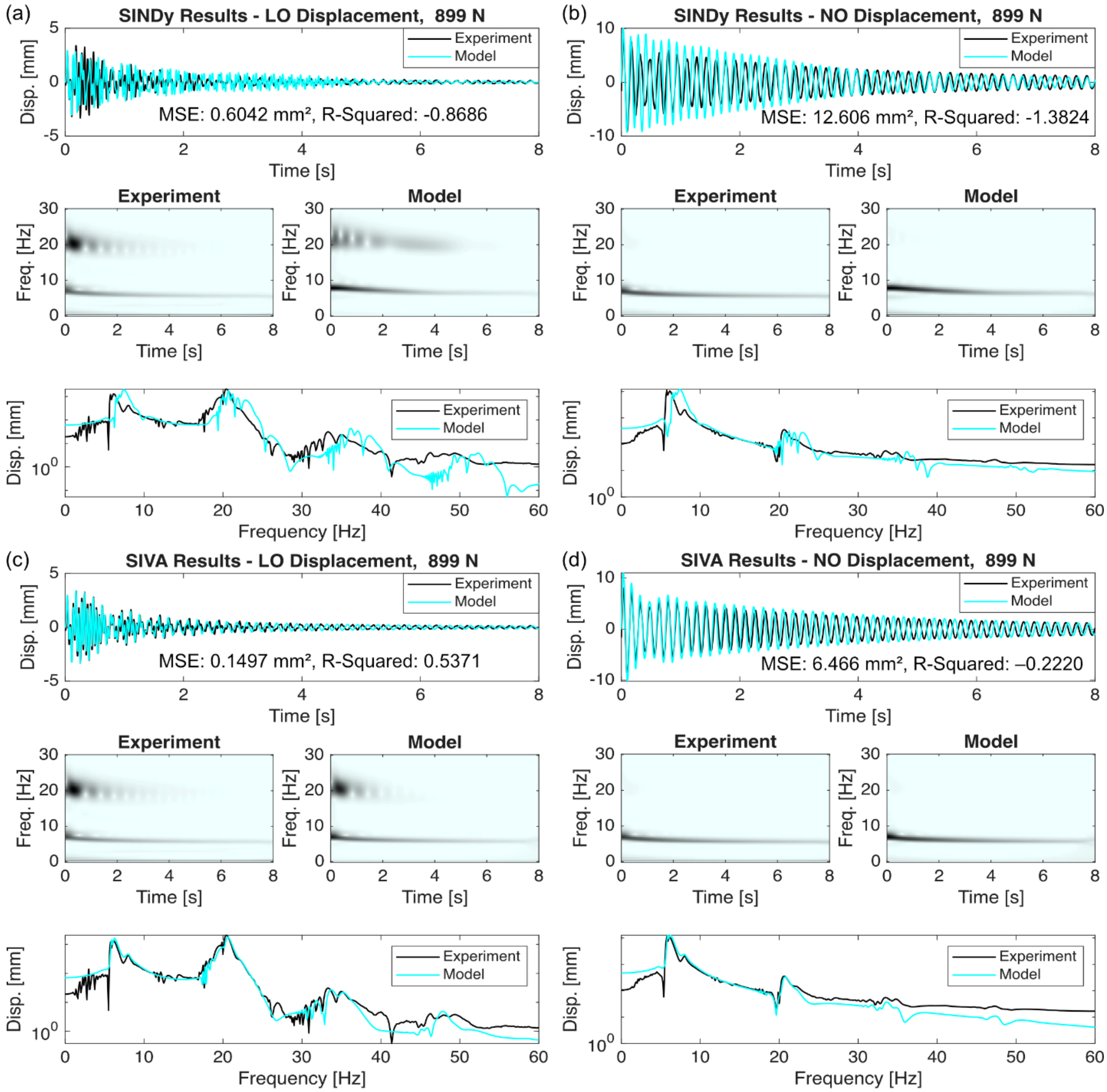

Figure 8. Comparison of the measured and predicted responses using the model produced by SINDy and SIVA for the measurement case used for the identification. Time-series, CWT spectra, and Fourier spectra shown for SINDy for (a) the LO and (b) the NO, and for SIVA for the (c) the LO and (d) the NO.

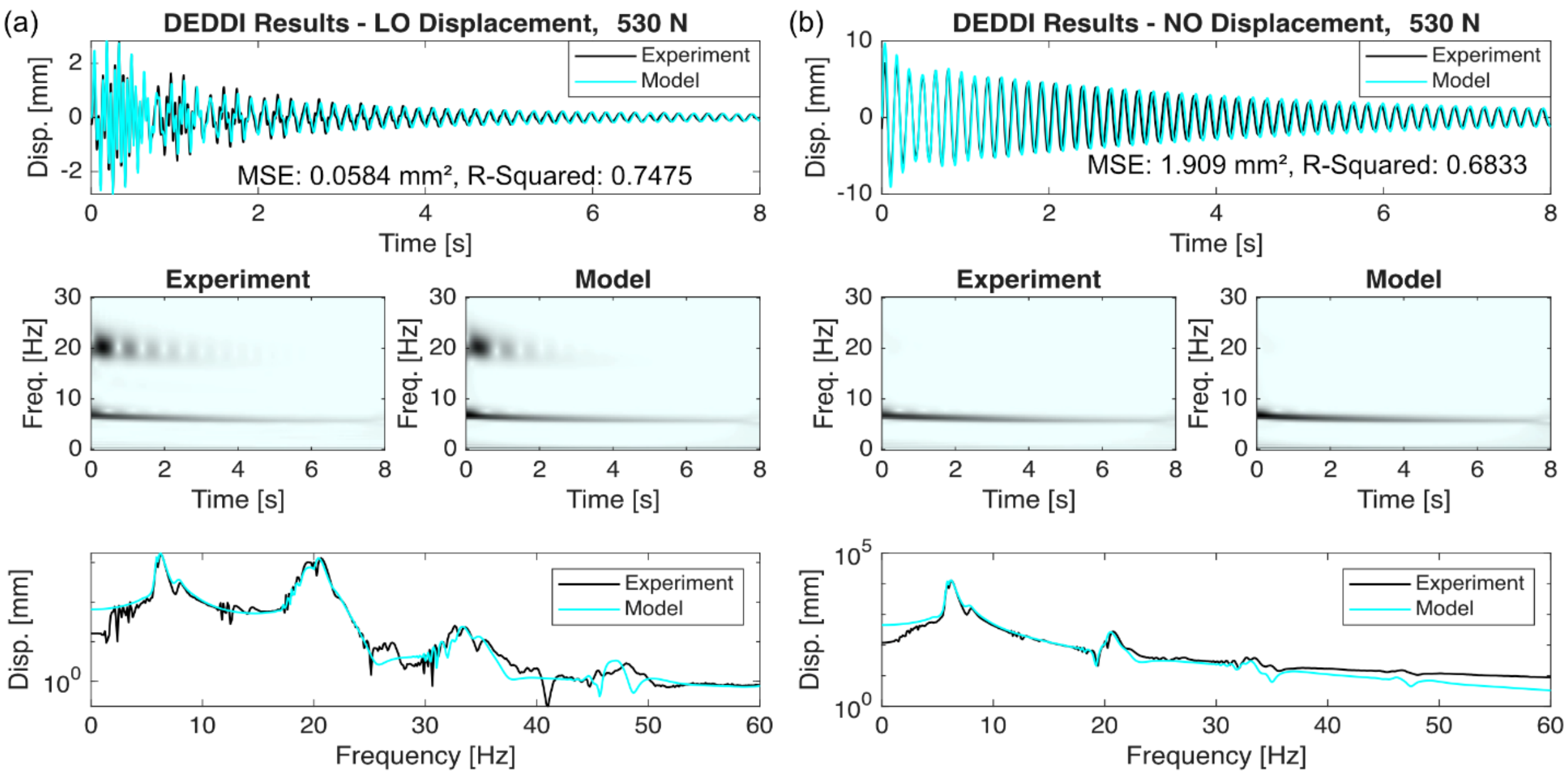

Figure 9. Comparison of the measured and predicted responses using the model produced by DEDDI for an excitation lower than the case used for identification. Time-series, CWT spectra, and Fourier spectra shown for (a) the LO and (b) the NO.

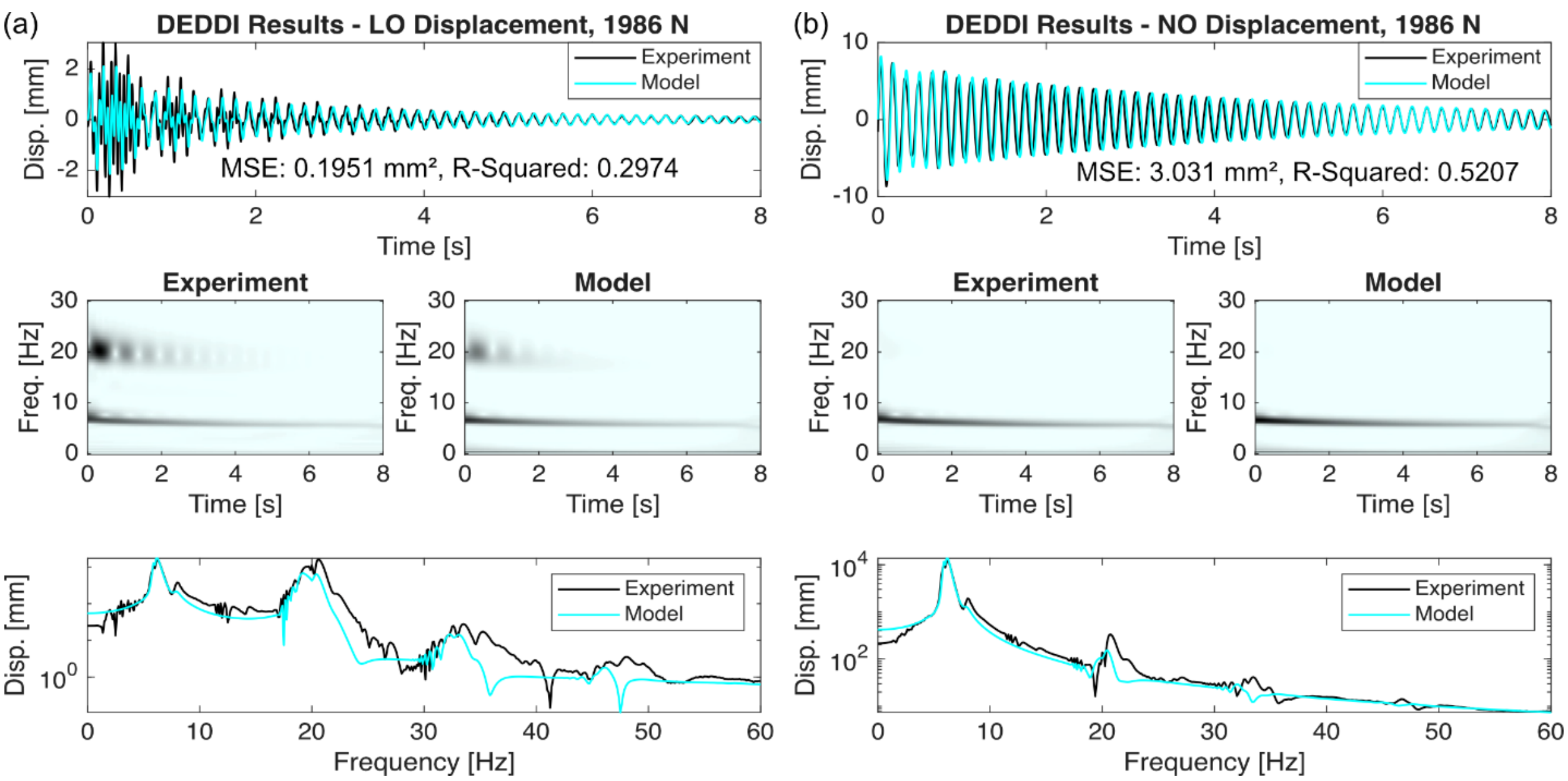

Figure 10. Comparison of the measured and predicted responses using the model produced by DEDDI for an excitation higher than the case used for identification. Time-series, CWT spectra, and Fourier spectra shown for (a) the LO and (b) the NO.

## 4. Concluding Remarks and Discussion

This paper introduced the Decomposition-based Energy-based Dual-phase Dynamics Identification (DEDDI) method for nonlinear system identification of MDOF systems. DEDDI combines WBEMD with the original EDDI framework to separately reconstruct the conservative and nonconservative internal forces acting on each DOF. First, WBEMD is applied to decompose the measured displacements into sets of IMFs. The IMFs are treated as SDOF oscillators and EDDI is applied to each individually. The resulting conservative and nonconservative forces are then summed across the IMFs to reconstruct the corresponding conservative and nonconservative internal forces acting on each DOF. These force components are subsequently used to identify independent mathematical models for stiffness and damping, yielding representative equations of motion for the system. The proposed framework was demonstrated using a two-story tower with nonlinear stiffness coupling the two floors. The results show that DEDDI reconstructs conservative and nonconservative internal forces and produces accurate mathematical models of the system dynamics. Comparison with SIVA and SINDy showed that DEDDI produced the most accurate model for the identification dataset, while SIVA exhibited slightly better generalization to the validation datasets.

While this work demonstrates the effectiveness of DEDDI, several theoretical questions and practical extensions remain open. First, there are two immediate extensions of DEDDI to be considered: forced response and continuous systems. For forced-response systems, it remains unclear how DEDDI should be applied when only a single harmonic is present in the response, particularly when the system is not synchronous at the forcing frequency. In the case of continuous systems, DEDDI could be applied to produce reduced-order models using discrete sets of sensors, while digital image correlation may enable full-field implementations at the expense of increased computational cost. An important open question is whether DEDDI can be extended beyond reduced-order models to identify continuum descriptions (i.e., PDEs) directly from measurements. Second, the theoretical basis for reconstructing the physical internal forces by summing independently identified IMF forces remains unknown. Although DEDDI accurately reconstructs the conservative and nonconservative internal forces for the systems considered here, the conditions under which this reconstruction is valid have not yet been established. In particular, computing the kinetic energy using decomposed velocities introduces cross terms that are neglected in the present formulation, yet the reconstructed forces remain physically meaningful. Explaining why these cross terms do not prevent accurate force reconstruction, and determining when they do, constitutes an important direction for future research. Lastly, further research is needed to understand the influence of internal resonances on DEDDI and to determine whether the associated energy transfers can be incorporated explicitly into the EDDI energy balance when applied to IMFs.

More broadly, DEDDI raises a fundamental question regarding decomposition-based system identification. Under what conditions can independently identified component forces be recombined to recover the governing force structure of the original physical system? Answering this question would establish a theoretical foundation for decomposition-based system identification methods and clarify the limits of their applicability.